\documentclass[a4paper, 12pt, twoside, reqno, openright]{amsart}
\usepackage{amsmath}
\usepackage{amssymb}
\usepackage{amsthm}
\usepackage{amscd}
\usepackage{amsopn}
\usepackage[german, french, UKenglish]{babel}
\usepackage{pifont}
\usepackage{lmodern}
\usepackage[cal=boondox, bb=ams, scr=boondox]{mathalfa}
\usepackage[top=1.5in, bottom=1.2in, left=1.2in, right=1.2in]{geometry}
\usepackage[all,cmtip]{xy}

\newtheoremstyle{Teorema}{5pt}{5pt}{\it}{}{\bf}{.}{ }{}
\theoremstyle{Teorema}
\newtheorem{Theorem}{Theorem}

\newtheorem*{Thrm}{Theorem}

\newtheorem*{ThrmA}{Theorem A}

\newtheorem*{Corol}{Corollary}
\newtheorem{Corollary}[Theorem]{Corollary}

\newtheorem{Definition}[Theorem]{Definition}
\newtheorem{Lemma}[Theorem]{Lemma}
\newtheoremstyle{Annotazione}{5pt}{5pt}{\rm}{}{\bf}{.}{ }{}
\theoremstyle{Annotazione}

\newtheorem*{Ackno}{Acknowledgements}

\newtheorem*{Outline}{Outline}

\def\Spec{\operatorname{Spec}}

\def\Aut{\operatorname{Aut}}
\def\Out{\operatorname{Out}}

\def\Autext{\operatorname{Out}}

\def\qbar{\overline{\bbq }}
\def\Ebar{\overline{E}}

\def\SL{\mathrm{SL}}

\def\GL{\mathrm{GL}}

\def\bbz{\mathbb{Z}}

\def\bbq{\mathbb{Q}}
\def\bbc{\mathbb{C}}

\def\bbp{\mathbb{P}}

\def\Gal{\operatorname{Gal}}

\makeatletter
\def\blfootnote{\xdef\@thefnmark{}\@footnotetext}
\makeatother

\begin{document}

\title{On copies of the absolute Galois group in $\Out\hat{F}_2$}
\author{Robert A. Kucharczyk}
\address{Universit\"{a}t Bonn\\ Mathematisches Institut\\ Endenicher Allee 60\\ 53115 Bonn\\ Germany}
\email{rak@math.uni-bonn.de}
\thanks{Research supported by the European Research Council}
\begin{abstract}
In this article we consider outer Galois actions on a free profinite group of rank two, induced by the \'{e}tale fundamental group of a projective line minus three points or of a pointed elliptic curve over a number field. Under mild technical assumptions their respective images uniquely determine the curves and the number fields.
\end{abstract}
\maketitle

\thispagestyle{empty}

\section{Introduction}

\noindent In 1979 G.V. Bely\u{\i} famously proved the following result, see \cite[Theorem~4]{MR534593}:
\begin{Thrm}[Bely\u{\i}]
Let $X$ be a smooth projective algebraic curve defined over $\qbar$. Then there exists a morphism $f\colon X\to\bbp^1_{\qbar }$ of algebraic curves over $\qbar$ which is unramified outside $\{ 0,1,\infty \}\subset\bbp^1$.
\end{Thrm}
Here we think of $\qbar$ as the set of all algebraic numbers in $\bbc$.

One of the most important consequences of this theorem is the existence of a continuous injective homomorphism
\begin{equation}\label{BelyiHomomorphismOut}
\varrho_{01\infty }\colon G_{\bbq }\hookrightarrow\Autext\hat{F}_2,
\end{equation}
where $G_{\bbq }=\Gal (\qbar /\bbq )$ is the absolute Galois group of the rational numbers, $\hat{F}_2$ is the profinite completion of a free group on two letters and $\Autext$ denotes the outer isomorphism group. The map (\ref{BelyiHomomorphismOut}) is obtained from the short exact sequence of \'{e}tale fundamental groups
\begin{equation}\label{SesForP1MoinsTroisPoints}
1\to\pi_1(\bbp_{\qbar }^1\smallsetminus\{ 0,1,\infty\} ,\ast )\to\pi_1(\bbp_{\bbq }^1\smallsetminus\{ 0,1,\infty \},\ast )\to G_{\bbq }\to 1
\end{equation}
in which the kernel can be identified with the profinite completion of
$$\pi_1^{\mathrm{top}}(\bbp^1(\bbc )\smallsetminus\{ 0,1,\infty \} ,\ast )\simeq F_2.$$
From Bely\u{\i}'s theorem it is easy to deduce:
\begin{Corol}[Bely\u{\i}]
The map (\ref{BelyiHomomorphismOut}) is injective.
\end{Corol}
Note that since $G_{\bbq }$ is compact (\ref{BelyiHomomorphismOut}) is therefore a homeomorphism onto its image.

The proof, which can be found in \cite[\S 4]{MR534593}, is based on the following observations: by Bely\u{\i}'s theorem every algebraic curve over $\qbar$ is birational to some finite \'{e}tale covering of $\bbp^1\smallsetminus\{ 0,1,\infty \}$, but such coverings correspond to conjugacy classes of open subgroups of $\hat{F}_2$. Since $G_{\bbq }$ operates faithfully on birationality classes of algebraic curves over $\qbar$, the corollary follows.

Choosing a base point $\ast$ defined over $\bbq$ we obtain a splitting of the sequence (\ref{SesForP1MoinsTroisPoints}) and hence a lift of (\ref{BelyiHomomorphismOut}) to an injection
\begin{equation}\label{BelyiHomomorphismAut}
G_{\bbq }\hookrightarrow\Aut\hat{F}_2;
\end{equation}
the most popular base point is the tangential base point $\ast =\overrightarrow{01}$ as defined in \cite{MR1012168}. Alexander Grothendieck urged his fellow mathematicians in \cite{MR1483107} to study the image of (\ref{BelyiHomomorphismOut}) or (\ref{BelyiHomomorphismAut}) with the hope of arriving at a purely combinatorial description of $G_{\bbq }$. He gave a candidate for the image, known today as the (profinite) Grothendieck--Teichm\"{u}ller group $\widehat{\mathrm{GT}}\subset\Aut\hat{F}_2$ (see \cite{MR1483118} for an overview). By construction $G_{\bbq }\hookrightarrow\widehat{\mathrm{GT}}$, but the other inclusion remains an open conjecture.

There are, however, still other embeddings $G_K\hookrightarrow\Autext\hat{F}_2$ for each number field $K\subset\bbc$. For each elliptic curve $E$ over $K$ we set $E^{\ast }=E\smallsetminus\{ 0\}$ and obtain a short exact sequence
$$1\to\pi_1(E_{\qbar }^{\ast })\to\pi_1(E^{\ast })\to G_K\to 1$$
analogous to (\ref{SesForP1MoinsTroisPoints}). Choosing a basis $\mathcal{B}$ of $H_1(E(\bbc ),\bbz )$ we construct an identification $\pi_1(E^{\ast }_{\qbar })\simeq\hat{F}_2$ below, and hence an injection (cf.\ Theorem~\ref{OuterGaloisInjective} below)
$$\varrho_E=\varrho_{E,\mathcal{B}}\colon G_K\hookrightarrow\Autext\hat{F}_2.$$
We will actually require this basis to be positive, i.e.\ positively oriented for the intersection pairing.

\begin{ThrmA}
For $j=1,2$ let $K_j\subset\bbc$ be a number field, $E_j$ an elliptic curve over $K_j$ and $\mathcal{B}_j$ a positive basis of $H_1(E_j(\bbc ),\bbz )$. Assume that
$$\varrho_{E_1,\mathcal{B}_1}(G_{K_1})=\varrho_{E_2,\mathcal{B}_2}(G_{K_2})$$
as subgroups of $\Out\hat{F}_2$. Then $K_1=K_2$ and there exists an isomorphism $E_1\simeq E_2$ over $K_1$ sending $\mathcal{B}_1$ to $\mathcal{B}_2$.
\end{ThrmA}
It is necessary to assume that the bases are positive: let $\tau$ denote complex conjugation, let $K$ be a non-real number field with $\tau (K)=K$ and let $E$ be an elliptic curve over $K$ with $E$ not isomorphic to $\tau (E)$. Then complex conjugation defines a real diffeomorphism $E(\bbc )\to\tau (E)(\bbc )$ sending each positive basis $\mathcal{B}$ of $H_1(E(\bbc ),\bbz )$ to a negative basis $\tau (\mathcal{B})$ of $H_1(\tau (E)(\bbc ),\bbz )$, and $\varrho_{E,\mathcal{B}}$ and $\varrho_{\tau (E),\tau (\mathcal{B})}$ have the same image.
\begin{Outline}
After recalling some preparatory material in section~2 we will prove Theorem~A in section~3 and finally draw some easy consequences in section~4.
\end{Outline}
\begin{Ackno}
The author wishes to thank Ursula Hamenst\"{a}dt for helpful comments on an earlier version of this article.
\end{Ackno}

\section{Some anabelian geometry}

\noindent We recall some facts about \'{e}tale fundamental groups of hyperbolic curves over number fields.
\begin{Definition}
Let $k$ be a field and $Y$ a smooth curve over $k$. Let $X$ be the smooth projective completion of $Y$ and $S=X(\overline{k})\smallsetminus Y(\overline{k})$; let $g$ be the genus of $X$ and $n$ the cardinality of $S$. Then $Y$ is called \emph{hyperbolic} if $\chi (Y)=2-2g-n<0$.
\end{Definition}
If $k\subseteq\bbc$ then $Y$ is hyperbolic if and only if the universal covering space of $Y(\bbc )$ is biholomorphic to the unit disk. Both $\bbp^1$ minus three points and an elliptic curve minus its origin are hyperbolic.

Now assume that $k=K\subset\bbc$ is a number field. By \cite[XIII~4.3]{MR0354651} the sequence
\begin{equation}\label{SesGeneralHyperbolicCurve}
1\to\pi_1(Y_{\qbar },\ast )\to\pi_1(Y,\ast )\to G_K\to 1
\end{equation}
induced by the ``fibration'' $Y_{\qbar }\to Y\to\Spec K$ is exact. By the usual group-theoretic constructions this sequence defines a homomorphism
\begin{equation}\label{OuterActionGeneralHyperbolic}
G_K\to\Out\pi_1(Y_{\qbar },\ast ),
\end{equation}
and the group $\pi_1(Y_{\qbar },\ast )$ is the profinite completion of $\pi_1^{\mathrm{top}}(Y(\bbc ),\ast )$, which is either a free group (in the affine case) or can be presented as
\begin{equation*}
\langle a_1,\ldots ,a_g,b_1,\ldots ,b_g\mid [a_1,b_1]\cdots [a_g,b_g]=1\rangle .
\end{equation*}
\begin{Theorem}\label{OuterGaloisInjective}
The homomorphism (\ref{OuterActionGeneralHyperbolic}) is injective.
\end{Theorem}
\begin{proof}
This is \cite[Theorem~C]{MR2895284}.
\end{proof}
We also note for later use that the sequence (\ref{SesGeneralHyperbolicCurve}) can be reconstructed from~(\ref{OuterActionGeneralHyperbolic}):
\begin{Lemma}\label{LemmaHomToOutDeterminesExtension}
Let $G$ be a profinite group, and let $\pi$ be a profinite group which is isomorphic to the \'{e}tale fundamental group of a hyperbolic curve over $\bbc$. Let $\varphi\colon G\to\Out\pi$ be a continuous group homomorphism. Then there exists a short exact sequence of profinite groups
$$1\to \pi\to H\to G\to 1$$
inducing $\varphi$, and it is unique in the sense that if another such sequence is given with $H'$ in the middle, then there exists an isomorphism $H'\to H$ such that the diagram
\begin{equation*}
\xymatrix{1\ar[r] & \pi \ar@{=}[d]\ar[r] & H\ar[d]^{\simeq }\ar[r] & G\ar@{=}[d] \ar[r] &1\\
1\ar[r] & \pi\ar[r] & H'\ar[r] & G \ar[r] &1
}
\end{equation*}
commutes.
\end{Lemma}
\begin{proof}
Let $\mathcal{Z}(\pi )$ denote the centre of $\pi$. The obstruction to the existence of such a sequence is a class in $H^3(G,\mathcal{Z}(\pi ))$ by \cite[Chapter~IV, Theorem~8.7]{MR0156879}, but since $\mathcal{Z}(\pi )$ is trivial by \cite[Proposition~18]{MR0354882}, the obstruction is automatically zero. Given the existence of one such sequence, the isomorphism classes of all such sequences are in bijection with $H^2(G,\mathcal{Z}(\pi ))=0$ by \cite[Chapter~IV, Theorem~8.8]{MR0156879}.
\end{proof}
We may safely ignore basepoints for the following reason: if $y,y'\in Y(\qbar )$ then there exists an isomorphism $\pi_1(Y_{\qbar },y)\simeq\pi_1(Y_{\qbar },y')$, canonical up to inner automorphisms. Hence the outer automophism groups of both are canonically identified. Furthermore, since both basepoints map to the same tautological base point of $\Spec k$, the whole sequence (\ref{SesGeneralHyperbolicCurve}) is changed only by inner automorphisms of the kernel when basepoints are changed within $Y(\qbar )$. So we drop basepoints from the notation in the sequel.

If $X$ and $Y$ are hyperbolic curves over a number field $K$ and $f\colon X\to Y$ is an isomorphism, we obtain a commutative diagram
\begin{equation}\label{IsomorphismOFSes}
\xymatrix{1\ar[r] & \pi_1(X_{\qbar }) \ar[d]^{\simeq }\ar[r] & \pi_1(X)\ar[d]^{\simeq }\ar[r] & G_K\ar@{=}[d] \ar[r] &1\\
1\ar[r] & \pi_1(Y_{\qbar })\ar[r] & \pi_1(Y)\ar[r] & G_K \ar[r] &1.
}
\end{equation}
\begin{Theorem}\label{TheoremTamagawa}
Let $K$ be a number field and $X$, $Y$ hyperbolic curves over $K$. Let $f\colon\pi_1(X)\to\pi_1(Y)$ be an isomorphism of fundamental groups commuting with the projections to $G_K$. Then $f$ is induced by a unique isomorphism of $K$-varieties $X\to Y$, and can be inserted into a commutative diagram of the form (\ref{IsomorphismOFSes}).
\end{Theorem}
\begin{proof}
This holds more generally for $K$ finitely generated over $\bbq$. It was conjectured by Grothendieck in \cite{MR1483108}, proved in the affine case by Tamagawa  in \cite[Theorem~0.3]{MR1478817} and in the projective case by Mochizuki in \cite{MR1720187}.
\end{proof}

\section{The Galois actions}

\noindent Let $F_2$ be the free group on two letters and let $\hat{F}_2$ be its profinite completion. Consider the following objects:
\begin{enumerate}
\item a number field $K\subset\bbc$,
\item an elliptic curve $E$ over $K$ and
\item a basis $\mathcal{B}$ of the homology group $H_1(E(\bbc ),\bbz )$.
\end{enumerate}
Let $E^{\ast }=E\smallsetminus\{ 0\}$, then $\pi_1^{\mathrm{top}}(E^{\ast }(\bbc))$ is a free group of rank two whose maximal abelian quotient can be identified with $H_1(E(\bbc ),\bbz)$. By the following lemma, the group isomorphism
\begin{equation*}
\bbz^2\to H_1(E(\bbc ),\bbz ),\quad (m,n)\mapsto mx+ny\text{ where }\mathcal{B}=(x,y)
\end{equation*}
 can be lifted uniquely to an outer isomorphism class
\begin{equation}\label{OuterBasisOfPiOneEStar}
F_2\dashrightarrow \pi_1^{\mathrm{top}}(E^{\ast }(\bbc )),
\end{equation} 
i.e.\ a group isomorphism which is well-defined up to inner automorphisms (which allows us to drop the basepoint for the fundamental group).
\begin{Lemma}\label{IsomorphismBetweenFreeAbelianInducesIsomorphismBetweenFree}
Let $F$ and $G$ be free groups of rank two, and let $f\colon F^{\mathrm{ab}}\to G^{\mathrm{ab}}$ be an isomorphism between their maximal abelian quotients. Then there exists an isomorphism $\tilde{f}\colon F\to G$ inducing $F$; it is uniquely determined by $f$ up to inner automorphisms of $F$.
\end{Lemma}
\begin{proof}
It is enough to prove this lemma in the case where $F=G=F_2$; but this is a reformulation of the well-known result that the natural map
\begin{equation*}
\Autext F_2\to\Aut (\bbz^2)=\GL (2,\bbz )
\end{equation*}
is an isomorphism.
\end{proof}
Since the profinite completion of $\pi_1^{\mathrm{top}}(E^{\ast }(\bbc ))$ can be identified with $\pi_1(E^{\ast }_{\qbar })$ we obtain an outer isomorphism class
\begin{equation}\label{OuterBasisOfPiOneGeomEStar}
\iota_{\mathcal{B}}\colon\hat{F}_2\dashrightarrow\pi_1(E^{\ast }_{\qbar }).
\end{equation}
Hence pulling back the Galois action on $\pi_1 (E^{\ast }_{\qbar })$ along (\ref{OuterBasisOfPiOneGeomEStar}) defines an injective homomorphism
\begin{equation}\label{RhoEDeclared}
\varrho_{E,\mathcal{B}}\colon G_K\hookrightarrow\Out\hat{F}_2.
\end{equation}
\begin{Lemma}\label{TopologicalAndGaloisIsomorphism}
Let $E$ be an elliptic curve over $\qbar$ and let $\sigma\in G_{\bbq }$. Let
$$f\colon\pi_1(E^{\ast })\dashrightarrow \pi_1(\sigma (E^{\ast }))$$
be an outer isomorphism class of profinite groups which can be obtained in each of the following ways:
\begin{enumerate}
\item it is the map of \'{e}tale fundamental groups induced via functoriality by the tautological isomorphism of schemes $t\colon E^{\ast }\to\sigma (E^{\ast })$;
\item it is the profinite completion of an outer isomorphism class
\begin{equation*}\pi_1^{\mathrm{top}}(E^{\ast }(\bbc ))\dashrightarrow\pi_1^{\mathrm{top}}(\sigma (E^{\ast })(\bbc ))
\end{equation*}
induced by an orientation-preserving isomorphism of real Lie groups 
$$h\colon E(\bbc )\to \sigma (E)(\bbc ).$$
\end{enumerate}
Then $\sigma$ is the identity, and so is the fundamental group isomorphism in \textup{(ii)}.
\end{Lemma}
\begin{proof}
Let $E[2]\subset E$ be the $2$-torsion subgroup and set $E^{\dagger }=E\smallsetminus E[2]$. The multiplication-by-$2$ map is a normal \'{e}tale covering $E^{\dagger }\to E^{\ast }$, therefore $\pi_1(E^{\dagger })$ is a normal open subgroup of $\pi_1(E^{\ast })$. From assumption (i) we see that $f$ maps $\pi_1(E^{\dagger })$ isomorphically to $\pi_1(\sigma (E^{\dagger }))$. Similarly $h$ maps $E[2]$ to $\sigma (E)[2]$, hence (i) and (ii) hold with every $\ast$ replaced by $\dagger$.

The quotient of $E^{\dagger }$ by the identification $x\sim -x$ is isomorphic over $\qbar$ to a scheme of the form $\bbp^1\smallsetminus\{ 0,1,\infty ,\lambda \}$, and we obtain a commutative diagram of schemes
\begin{equation}\label{OddDiagramSchemes}
\xymatrix{
E^{\dagger }\ar[r]^-{\wp }\ar[d]_t & \bbp^1_{\qbar }\smallsetminus\{ 0,1,\infty ,\lambda \} \ar[r]^-{\iota } &\bbp^1_{\qbar }\smallsetminus\{ 0,1,\infty \} \ar[d]_t\\
\sigma (E^{\dagger })\ar[r]_-{\wp } & \bbp^1_{\qbar }\smallsetminus\{ 0,1,\infty ,\sigma (\lambda )\} \ar[r]_-{\iota } &\bbp^1_{\qbar }\smallsetminus\{ 0,1,\infty \}\\
}
\end{equation}
where the horizontal maps are morphisms of $\qbar$-schemes and all vertical maps are base change morphisms along $\sigma\colon\qbar\to\qbar$. (The maps $\wp$ are not necessarily Weierstra\ss\ $\wp$-functions, but up to M\"{o}bius transformations on $\bbp^1$ they are, whence our notation.)

There is a very similar commutative diagram of topological spaces:
\begin{equation}\label{OddDiagramTopol}
\xymatrix{
E^{\dagger }(\bbc )\ar[r]^-{\wp }\ar[d]_h & \bbp^1(\bbc )\smallsetminus\{ 0,1,\infty ,\lambda \} \ar[r]^-{\iota } &\bbp^1(\bbc )\smallsetminus\{ 0,1,\infty \} \ar[d]_H\\
\sigma (E^{\dagger })(\bbc )\ar[r]_-{\wp } & \bbp^1(\bbc )\smallsetminus\{ 0,1,\infty ,\sigma (\lambda )\} \ar[r]_-{\iota } &\bbp^1(\bbc )\smallsetminus\{ 0,1,\infty \} \\
}
\end{equation}
where the horizontal maps are holomorphic maps between Riemann surfaces and the vertical maps are orientation-preserving homeomorphisms.

We claim that the horizontal compositions in (\ref{OddDiagramSchemes}) and (\ref{OddDiagramTopol}), i.e.\ the composite maps of the form
\begin{equation}\label{CompositionOfSchemeMorph}
E^{\dagger }\overset{\wp }{\to }\bbp^1\smallsetminus\{ 0,1,\infty ,\lambda \}\overset{\iota }{\to }\bbp^1\smallsetminus\{ 0,1,\infty \} ,
\end{equation}
induce surjections on fundamental groups. This can be checked in the topological case; the maps induced by (\ref{CompositionOfSchemeMorph}) on topological fundamental groups have the form
\begin{equation}\label{CompositionOfGroupMorph}
\operatorname{ker}s\to\langle a_1,a_2,a_3,a_4\mid a_1a_2a_3a_4=1\rangle \to \langle a_1,a_2,a_3\mid a_1a_2a_3=1\rangle
\end{equation}
where the second map is given by $a_4\mapsto 1$, and
$$s\colon\langle a_1,a_2,a_3,a_4\mid a_1a_2a_3a_4=1\rangle \to\bbz /2\bbz,\quad a_i\mapsto 1\bmod 2.$$
Hence $a_ja_4\mapsto a_j$ for $j=1,2,3$ under (\ref{CompositionOfGroupMorph}), and therefore (\ref{CompositionOfGroupMorph}) is surjective.

From this we deduce that the two diagrams (\ref{OddDiagramSchemes}) and (\ref{OddDiagramTopol}) induce the same commutative diagrams of outer homomorphism classes between the \'{e}tale fundamental groups: the groups are clearly the same, and so are the homomorphisms induced by the horizontal maps and by the vertical maps on the left. But since the maps $\iota\circ\wp$ induce surjections on fundamental groups, the vertical maps on the right also have to induce the same homomorphisms.

In particular the base change map $t$ induced by $\sigma\colon\qbar\to\qbar$ and the orientation-preserving homeomorphism $H$ define the same element in $\Autext\pi_1(\bbp^1_{\qbar }\smallsetminus \{ 0,1,\infty \} )$. But $H$ is homotopic to the identity, hence this element has to be trivial; and by Theorem~\ref{OuterGaloisInjective} for $Y=\bbp^1\smallsetminus\{ 0,1,\infty\}$ the automorphism $\sigma$ has to be trivial, too.
\end{proof}

We note a result closely related to Lemma~\ref{TopologicalAndGaloisIsomorphism}, see \cite{MR1781929}:

\begin{Theorem}[Matsumoto--Tamagawa]\label{TheoremMatsumotoTamagawa}
Let $E$ be an elliptic curve defined over a number field $K\subset\bbc$. Then the images of the outer Galois representation
$$\Gal (\qbar /K )\to\Autext\pi_1 (E^{\ast }_{\qbar })$$
and the profinite closure of the topological monodromy
$$\widehat{\SL (2,\bbz )}\to \Autext\pi_1(E^{\ast }_{\bbc })=\Autext\pi_1 (E^{\ast }_{\qbar })$$
intersect trivially.
\end{Theorem}

We also need one more result on isomorphisms, this time between Galois groups. Let $K$ and $L$ be number fields in $\bbc$, and assume that $\sigma\in G_{\bbq }$ satisfies $\sigma (K)=L$. Then we can define a group isomorphism
$$\Phi_{\sigma }\colon G_K\to G_L,\quad \tau\mapsto\sigma\tau\sigma^{-1}.$$

\begin{Theorem}[Neukirch--Uchida]\label{TheoremNeukirchUchidaPop}
Let $K,L\subset\bbc$ be number fields and let $\Phi\colon G_K\to G_L$ be a continuous group isomorphism. Then there exists a unique $\sigma\in G_{\bbq }$ with $\sigma (K)=L$ and $\Phi =\Phi_{\sigma }$.
\end{Theorem}
For the proof see \cite{MR547650}.

\begin{proof}[Proof of Theorem~A]
The bases $\mathcal{B}_j$ of $H_1(E_j(\bbc ),\bbz )$ define an orientation-pre\-serv\-ing isomorphism between these two cohomology groups, hence an orientation-pre\-serv\-ing isomorphism of \emph{real} Lie groups $h\colon E_1(\bbc )\to E_2(\bbc )$ and an isomorphism of profinite fundamental groups
$$h_{\ast }=\iota_{\mathcal{B}_2}^{-1}\circ\iota_{\mathcal{B}_1}\colon \pi_1(E_{1,\qbar }^{\ast })\to\pi_1(E_{2,\qbar }^{\ast }).$$
Since the representations $\varrho_{E_j}$ are injective there is a unique isomorphism of profinite groups $\Phi\colon G_{K_1}\to G_{K_2}$  such that $\varrho_{E_1}=\varrho_{E_2}\circ\Phi$. By Theorem~\ref{TheoremNeukirchUchidaPop} this has to be of the form $\Phi_{\sigma }$ for a unique isomorphism $\sigma\in G_{\bbq }$ with $\sigma (K_1)=K_2$. We shall construct an isomorphism $\sigma (E_1)\to E_2$ of elliptic curves over $K_2$.

Consider the short exact homotopy sequences for the three varieties $\sigma (E_1)$, $E_1$, $E_2$ over their respective base fields; they can be completed to the following commutative diagram:
\begin{equation*}
\xymatrix{
1 \ar[r] & \pi_1(\sigma (\Ebar^{\ast }_1)) \ar[r]\ar[d]^{\simeq }_{m_{\ast }} & \pi_1(\sigma (E_1^{\ast })) \ar[r]\ar[d]^{\simeq }_{m_{\ast }} & G_{K_2} \ar[r]\ar[d]_{\Phi_{\sigma }^{-1}} & 1\\
1 \ar[r] & \pi_1(\Ebar^{\ast }_1) \ar[r]\ar[d]^{\simeq }_{h_{\ast }} & \pi_1(E_1^{\ast }) \ar[r] & G_{K_1} \ar[r]\ar[d]_{\Phi_{\sigma }} & 1\\
1 \ar[r] & \pi_1(\Ebar^{\ast }_2) \ar[r] & \pi_1(E_2^{\ast }) \ar[r] & G_{K_2} \ar[r] & 1.
}
\end{equation*}
Here the lower rectangle commutes trivially by exactness of the rows, and the upper two squares commute by functoriality of the fundamental group. From Lemma~\ref{LemmaHomToOutDeterminesExtension} we obtain an isomorphism $F\colon\pi_1(E_1^{\ast })\to\pi_1(E_2^{\ast })$ that makes the resulting diagram commute:
\begin{equation*}
\xymatrix{
1 \ar[r] & \pi_1(\sigma (E^{\ast }_{1,\qbar })) \ar[r]\ar[d]^{\simeq }_{m_{\ast }} & \pi_1(\sigma (E_1^{\ast })) \ar[r]\ar[d]^{\simeq }_{m_{\ast }} & G_{K_2} \ar[r]\ar[d]_{\Phi_{\sigma }^{-1}} & 1\\
1 \ar[r] & \pi_1(E^{\ast }_{1,\qbar }) \ar[r]\ar[d]^{\simeq }_{h_{\ast }}& \pi_1(E_1^{\ast }) \ar[r]\ar[d]_F & G_{K_1} \ar[r]\ar[d]_{\Phi_{\sigma }} & 1\\
1 \ar[r] & \pi_1(E^{\ast }_{2,\qbar }) \ar[r] & \pi_1(E_2^{\ast }) \ar[r] & G_{K_2} \ar[r] & 1.
}
\end{equation*}
By Theorem~\ref{TheoremTamagawa} the group isomorphism $F\circ m_{\ast }\colon\pi_1(\sigma (E_1^{\ast }))\to\pi_1(E_2^{\ast })$ must be induced by a unique isomorphism $g\colon\sigma (E_1)\to E_2$ of $K_2$-schemes. But this means that
\begin{equation*}
m_{\ast }\colon\pi_1(\sigma (E_{1,\qbar }^{\ast }))\to\pi_1(E_{1,\qbar }^{\ast })
\end{equation*}
is induced by the orientation-preserving homeomorphism
\begin{equation*}
h^{-1}\circ g^{\mathrm{an}};
\end{equation*}
by Lemma~\ref{TopologicalAndGaloisIsomorphism} we find that $\sigma$ must be the identity, so $K_1=K_2$ and $g$ is the desired isomorphism.
\end{proof}

\section{Concluding remarks}

\noindent From Theorem~A we can easily deduce several analogous statements. Recall that two subgroups $H',H''$ of a group $G$ are called directly commensurable if $H'\cap H''$ has finite index both in $H'$ and in $H''$; they are called widely commensurable if $gH'g^{-1}$ and $H''$ are directly commensurable for some $g\in G$.
\begin{Corollary}
For $j=1,2$ let $K_j\subset\bbc$ be a number field, $E_j$ an elliptic curve over $K_j$ and $\mathcal{B}_j$ a positive basis of $H_1(E_j(\bbc ),\bbz )$. Let $I_j$ be the image of $\varrho_{E_j,\mathcal{B}_j}\colon G_{K_j}\to\Out\hat{F}_2$.
\begin{enumerate}
\item $I_1=I_2$ if and only if $K_1=K_2$ and there exists an isomorphism $E_1\simeq E_2$ over $K_1$ sending $\mathcal{B}_1$ to $\mathcal{B}_2$.
\item $I_1$ and $I_2$ are conjugate in $\Out\hat{F}_2$ if and only if there exists a field isomorphism $\sigma\colon K_1\to K_2$ such that $\sigma (E_1)\simeq E_2$ as elliptic curves over $K_2$.
\item $I_1$ and $I_2$ are directly commensurable if and only if there exists an isomorphism $E_{1,\bbc }\to E_{2,\bbc }$ sending $\mathcal{B}_1$ to $\mathcal{B}_2$.
\item $I_1$ and $I_2$ are widely commensurable if and only if $j(E_1)$ and $j(E_2)$ lie in the same $G_{\bbq }$-orbit in $\qbar$, where $j$ denotes the $j$-invariant from the classical theory of elliptic curves.
\end{enumerate}
\end{Corollary}
\begin{proof}
(i) is Theorem~A and (ii) follows easily from Theorems \ref{TheoremTamagawa} and~\ref{TheoremNeukirchUchidaPop}. For (iii) we can find an open subgroup $G_{L_j}$ of each $G_{K_j}$ such that these two subgroups have the same image; we can then apply (i) to $E_i\otimes_{K_i}L_i$. Vice versa any isomorphism between two elliptic curves over $\bbc$ that admit models over number fields must already be defined over some number field. (iv) follows similarly from (ii).
\end{proof}
\begin{Corollary}
Let $K$ be a number field, $E$ an elliptic curve over $K$ and $\mathcal{B}$ a basis of $H_1(E_j(\bbc ),\bbz )$. Then $\varrho_{E,\mathcal{B}}(G_K)$ and $\varrho_{01\infty }(G_{\bbq })$ are not widely commensurable in $\Out\hat{F}_2$.
\end{Corollary}
\begin{proof}
Assume they were widely commensurable; after enlarging the fields of definition $K$ and $\bbq$ to some suitable number fields $L_1$, $L_2$ the two Galois images would actually be conjugate in $\Out\hat{F}_2$. As in the proof of Theorem~A we would obtain an isomorphism $\sigma\in G_{\bbq }$ with $\sigma (L_1)=L_2$ and a commutative diagram
\begin{equation*}
\xymatrix{
1\ar[r] & \pi_1(\sigma (E_{\qbar }^{\ast }))\ar[d]^{\simeq }\ar[r] & \pi_1(\sigma (E)_{L_2})\ar[d]^{\simeq }\ar[r] & G_{L_2} \ar@{=}[d]\ar[r] & 1\\
1\ar[r] & \pi_1(\bbp^1_{\qbar }\smallsetminus\{ 0,1,\infty \})\ar[r] & \pi_1(\bbp^1_{L_2}\smallsetminus\{ 0,1,\infty \})\ar[r] & G_{L_2} \ar[r] & 1,
}
\end{equation*}
hence by Theorem~\ref{TheoremTamagawa} an isomorphism $\sigma (E^{\ast }_{\qbar })\to \bbp^1_{\qbar }\smallsetminus\{ 0,1,\infty \}$, which is absurd.
\end{proof}

\providecommand{\bysame}{\leavevmode\hbox to3em{\hrulefill}\thinspace}
\providecommand{\MR}{\relax\ifhmode\unskip\space\fi MR }
\providecommand{\MRhref}[2]{%
  \href{http://www.ams.org/mathscinet-getitem?mr=#1}{#2}
}
\providecommand{\href}[2]{#2}

\end{document}